\documentclass[12pt]{amsart}
\usepackage{enumerate, amsmath, amssymb, hyperref}

\renewcommand{\phi}{\varphi}

\DeclareMathOperator{\dist}{dist}

\DeclareMathOperator{\sgn}{sgn}
\DeclareMathOperator*{\osc}{osc} 
\newtheorem{theorem}{Theorem}[section]
\newtheorem{lemma}[theorem]{Lemma}
\newtheorem{prop}[theorem]{Proposition}
\newtheorem{cor}[theorem]{Corollary}
\theoremstyle{definition}
\newtheorem{definition}[theorem]{Definition}
\newtheorem{example}[theorem]{Example}

\numberwithin{equation}{section}
\theoremstyle{remark}
\newtheorem{remark}[theorem]{Remark}

\title[Global invertibility of Sobolev mappings]{Global invertibility of Sobolev mappings with prescribed homeomorphic boundary values}
 \author[S. Traver]{Sabrina Traver}
 \address{Department of Mathematics, Syracuse University, Syracuse,
 NY 13244}
 \email{smtraver@syr.edu}
 \thanks{I would like to thank my advisor, Jani Onninen, for all of his invaluable discussions, guidance, and support throughout the duration of this project, as well as for the partial financial support provided through his NSF grant DMS-2453853. I would also like to thank Leonid Kovalev for his helpful comments, especially with regard to the degree function.}



\begin{document} 
\begin{abstract} 
Let $X, Y \subset \mathbb{R}^n$ be Lipschitz domains, and suppose there is a homeomorphism $\varphi \colon \overline{X} \to \overline{Y}$. We consider the class of Sobolev mappings $f \in W^{1,n} (X, \mathbb{R}^n)$ with a strictly positive Jacobian determinant almost everywhere, whose Sobolev trace coincides with $\varphi$ on $\partial X$. We prove that every mapping in this class extends continuously to $\overline{X}$ and is a monotone (continuous) surjection from $\overline{X}$ onto $\overline{Y}$ in the sense of C.~B.~Morrey. As monotone mappings, they may squeeze but not fold the reference configuration $X$. This behavior reflects weak interpenetration of matter, as opposed to folding, which corresponds to strong interpenetration. Despite allowing weak interpenetration of matter, these maps are globally invertible, generalizing the pioneering work of J.M. Ball. \end{abstract}
\maketitle

\section{Introduction}\label{sec1}
Throughout this paper, $X$ and $Y$ will represent bounded Lipschitz domains in $\mathbb R^n$.

Geometric Function Theory, particularly within variational frameworks~\cite{Astala, Hencl, Tadeuszbook, TMO}, and the theory of Nonlinear Elasticity~\cite{Ciarletbook, MardsenHughes, Tadeuszbook} studies elastic bodies $X$ and their deformations $f: X \to \mathbb R^n$. Properties of these deformations correspond to physical phenomena, making some classes of mappings better suited to the study of materials science. Most notably, well-behaved deformations are injective (or at least injective a.e.) because injectivity is essential for ensuring non-interpenetration of matter.

However, even $C^1$-mappings $f$ with strictly positive Jacobian determinant $J_f(x) = \det Df(x) >0$, which are locally injective by the classical inverse function theorem, may fail to be injective on the whole domain. This failure may arise, for example, when different ends of the domain are mapped to overlapping subsets of the target. The result is self-interpenetration of matter, but it can be shown using Degree Theory that when the $C^1$ map $f$ coincides with an orientation-preserving homeomorphism $\varphi: \overline{X} \to \overline{Y}$ on the boundary of the domain $X$, this cannot occur. More precisely, $C^1$ mappings $f: \overline{X} \to \mathbb R^n$ with strictly positive Jacobian determinant are homeomorphisms from $\overline{X}$ onto $\overline{Y}$ when $f|_{\partial X} = \varphi|_{\partial X}$.

Replacing the $C^1$ maps, $f$ and $\varphi$, with Sobolev maps, injectivity is generally lost (see example~\ref{mtexample}), but a weaker version of injectivity (injectivity a.e.) is still maintained. Moreover, mappings that are injective a.e., also known as globally invertible, are sufficient for physical applications~\cite{ball, Ciarlet, Kromer}. We say a map $f: X \to Y$ is \emph{globally invertible} if \begin{equation}\label{globalinvertibility} |\{y \in Y: f^{-1}(y) \text{ is a single point}\}| = |Y| 
\end{equation}

Equation (\ref{globalinvertibility}) is often referred to as global invertibility~\cite{ball} or invertibility for a.e. $y \in \overline{Y}$ because when (\ref{globalinvertibility}) holds, then a.e. $y \in \overline{Y}$ has a unique preimage. 

Each map in the following class of deformations 
\[\mathcal{A}^p_{\varphi} = \{f \in W^{1,p}(X, \mathbb R^n) : J_f > 0  \text{ a.e.}, f-\varphi \in W^{1,p}_0(X, \mathbb R^n)\}.\] is globally invertible for $p \ge n$. The $p > n$ case can be found in~\cite{ball}, and, as we will show in Section~\ref{app} the $p=n$ case follows from Theorem~\ref{thm-main}. When $p=n$ such maps are also monotone in the sense of Morrey~\cite{morrey}:
\begin{definition}\label{monotone}
    Given compact, connected metric spaces, $X$ and $Y$, a continuous map $f:X \to Y$ is \emph{monotone} if $f^{-1}(y_0) $ is connected for all $y_0 \in Y$.
\end{definition} 
Monotonicity allows for \emph{weak interpenetration of matter}. That is, roughly speaking, squeezing a portion of the material to a point can occur but not folding or tearing.
We now introduce our main result:
\begin{theorem}\label{thm-main}
    Let $X$ and Y be bounded Lipschitz domains in $\mathbb R^n$ and $\varphi:\overline{X} \to \overline{Y}$ be a given homeomorphism in $W^{1,n}(X, \mathbb R^n)$. If $f \in W^{1,n}(X, \mathbb R^n)$ with $J_f > 0$ a.e. such that $f- \varphi \in W^{1,n}_0(X, Y) $ then 
    \begin{enumerate}[(1)]
        \item $f$ has a continuous extension to the boundary of $X$, which we still denote $f$
        \item $f$ maps $\overline{X}$ onto $\overline{Y}$
         \item $f: \overline{X} \to \overline{Y}$ is a monotone map
         \item $f: \overline{X} \to \overline{Y}$ is globally invertible
    \end{enumerate}
\end{theorem}
\begin{cor}\label{cor-injective}
    Let $f$ be as in Theorem~\ref{thm-main}. Then there is a set $X_f \subset X$ of full measure such that $f$, when restricted to $X_f$, is injective.
\end{cor}

The $W^{1,p}$ version of Theorem~\ref{thm-main} for $p > n$ was proven by Ball in 1981~\cite{ball}. Though the results are similar, the conformally invariant case $p=n$ has one significant difference: for $p>n$, the boundary regularity follows directly from continuity estimates for arbitrary Sobolev maps. In the critical case $p=n$, such estimates break down, so we devise new estimates by analyzing oscillation and applying a Sobolev embedding theorem. 

Even in the $p=n$ case, slight changes to the class of maps may produce different results. As an example, Ball, Iwaniec, and Onninen~\cite{balliwanieconninen} previously studied weakly convergent sequences of $W^{1,n}-$Sobolev homeomorphisms in the traction free setting. When boundary data is not fixed, they showed that uniform equicontinuity may fail without additional topological or variational assumptions~\cite{balliwanieconninen}. In contrast, the present paper relies on homeomorphic boundary data to prove many properties of the map in Theorem~\ref{thm-main}, such as continuity up to the boundary and global invertibility.

In section~\ref{app}, we apply our result to solve an energy minimization problem. Although the classical results of Ioffe~\cite{ioffe} could potentially be applied to this and similar settings, the singular behavior of the integrand in combination with the assumption $J_f >0 $ a.e. means a direct proof is better suited to demonstrate the unique features of this functional. Indeed, the proof relies on convexity, weak $L^1_{\text{loc}}$ convergence of Jacobians, and an application of the Monotone Convergence Theorem.

The last section,~\ref{examples}, is reserved for examples. In particular, we show the strict inequality $J_f > 0$ a.e. is essential. In Example~\ref{ex-ctyonto}, if we allow $J_f \ge 0$, Theorem~\ref{thm-main}(1) and (3) may fail, and even if $f$ is also continuous up to the boundary, as in Example~\ref{ex-ctybdry}, Theorem~\ref{thm-main}(2) may fail. 

\section{Preliminaries}
The domains in this preliminaries section need not be Lipschitz. We say $X$ is a Lipschitz domain if for each $x \in \partial X$ there exists a neighborhood $U_x$ of $x$ and a Cartesian coordinate system $\hat{x} = \{x_1, \dots , x_n\}$ in $U_x$ such that $X \cap U_x = \{\hat{x} \in U_x : x_n > f(x_1, \dots , x_{n-1})\}$ for some Lipschitz continuous function $f: \mathbb R^{n-1} \to \mathbb R$.  
  
 We begin by introducing the Lusin $N$ condition, which, when considering physical applications, ensures that our deformation $f$ cannot create material from nothing. Mathematically, for instance, its value is in its ability to relate the Jacobian of our map $f$ to a degree function, which counts the number of preimages under $f$ of a given point, while accounting for orientation.
\begin{definition}
   Let $X \subset \mathbb R^n$ be open. A map $f:X \to \mathbb R^n$ satisfies the \emph{Lusin $N$ condition} if for every $E \subset X$ with $|E| =0$, $|f(E)|=0$
 \end{definition}

If the preimage of a null set is a null set, we say $f$ has the Lusin $N^{-1}$ condition.

\begin{remark}\label{conditionN}
 Any $f \in W^{1,n}(X, \mathbb R^n)$ with $J_f > 0$ a.e. has both the Lusin $N$ and $N^{-1}$ conditions. The former can be seen by applying~\cite[Theorem 4.5]{Hencl} with $|E|=0$ and the latter is proven in~\cite[Theorem 5.32]{Fonseca}.  
\end{remark}

To extend the map $f$ to $\partial X$ continuously, we prove that each component of the extended map is weakly monotone as defined below. We denote $(f-g)^+ := \max\{f-g, 0\}$.
\begin{definition}
    A real valued function $u \in W^{1,1}(X)$ is \emph{weakly monotone} if for every ball $B \subset X$ and all constants $m \le M$ such that $(u-M)^+ - (m-u)^+ \in W^{1,1}_0(B)$ implies that $m \le u(x) \le M$ for almost every $x \in B$. 
\end{definition}

One of the central tools used to prove properties of our extended map $f: \overline{X} \to \overline{Y}$ is the topological degree

 \begin{definition}\label{degree}
    Let $X$ be a bounded, open subset of $ \mathbb R^n$ and let $f: X \to \mathbb R^n$ be a continuous map. Then for any subset $B$ compactly contained in $X$ and $p \in \mathbb R^n \setminus (f(\partial B))$, the \emph{degree of $f$ at $p$ with respect to $B$} is 
    \[\deg (f, p, B) := \sum_{x \in \psi^{-1}(p) \cap B} \sgn J_{\psi}(x)\]
    for any $C^1$ map $\psi: \overline{B} \to \mathbb R^n$ such that $|\psi(x) - f(x)| < \dist(p, f(\partial B))$ for all $x \in \overline{X}$. 
\end{definition}
See~\cite{Fonseca, Schwartz} for properties of this degree function. The Lusin $N$ condition pairs nicely with the topological degree as seen in the following theorem
\begin{theorem}[\cite{Fonseca}, Theorem 5.35]\label{Fonseca}
    Let $D \subset \mathbb R^n$ be a bounded, open set, and let $f \in W^{1,n}(D, \mathbb R^n)$ be a mapping such that $J_f>0$ a.e. $x \in D$. If $v \in L^{\infty}(\mathbb R^n)$ then for every open set $G \subset D$ such that $|\partial G|=0$, we have 
\[\int_G v \circ f(x) J_f(x) dx = \int_{\mathbb R^n} v(y) \deg(f, y, G) dy\]
\end{theorem}
We modify this equality in section~\ref{sec-properties} to help prove various properties of our extended map $f$.

\section{Topological lemmas}\label{toplemmas}
Before beginning the proof of the main result, we present some lemmas that are purely topological. We denote an $n$ dimensional ball (and $(n-1)$-dimensional sphere in $\mathbb R^n$) with center $x$ and radius $r$ as $B_r(x)$ (and $S^{n-1}_r(x)$ respectively). 

\begin{lemma}\label{lem-component}
    Let $f:X \to \mathbb R^n$ be continuous, $y \in \mathbb R^n$, and let $K \subset X$ be a component of $f^{-1}(y)$ such that $\overline{K} \subset X$. Then 
    \begin{enumerate}
        \item there exists $r_K >0$ such that for every $r \in (0, r_K)$, the component $U_r$ of $f^{-1}(B_r(y))$ containing $K$ satisfies $\overline{U_r} \subset X$.
        \item if $K' \neq K$ is another component of $f^{-1}(y)$ such that $\overline{K'} \subset X$ then there exists $r_{K,K'} >0$ such that for every $r \in (0, r_{K,K'})$ the sets $K$ and $K'$ are contained in different components of $f^{-1}(B_r(y))$
    \end{enumerate}
\end{lemma}
\begin{proof} (1)
    Since $K$ is a component of $f^{-1}(y)$, it is contained in a single component of $f^{-1}(B_r(y))$. Suppose for the sake of contradiction that we can build a sequence $r_1 > r_2 > \cdots $ such that $r_i \to 0$ and $\overline{U_{r_i}} \cap \partial X \neq \emptyset$ for all $i$. 
    
    Let $C:= \cap_i \overline{U_{r_i}} $ and observe that $K \subset C$ and $C \cap \partial X \neq \emptyset$. Moreover, $\overline{U_{r_i}}$ is a decreasing sequence of compact, connected sets, so its intersection is nonempty, compact, and connected.

    As $f$ is continuous, $f(\overline{U_{r_i}} \cap X) \to f(C \cap X)$ as $r_i \to 0$, and since $f(\overline{U_{r_i}} \cap X) \subset \overline{B_{r_i}(y)}$, it follows that $f(C \cap X) \subset \{y\}$.

Now, $K$ is a connected component of $f^{-1}(y)$, and it follows that the connected component of $C \cap X$ containing $K$ must equal $K$. Then we can find an open neighborhood $U \supset K$ such that $\overline{U} \cap ((C \cap X)\setminus K) = \emptyset$, and since $\overline{K} \subset X$, there exists an open neighborhood $V$ of $K$ such that $\overline{V} \subset X$. 
    
    Now, $U\cap V$ and $C\setminus (\overline{U} \cup \overline{V})$ forms a separation of $C$, but $C$ is connected, so this is a contradiction. 

    (2) Suppose there exists a sequence $r_1 > r_2 > \cdots $ such that $r_i \to 0$ and components $U_i \subset f^{-1}(B_{r_i}(y))$ such that $U_i \cap K \neq \emptyset$ and $U_i \cap K' \neq \emptyset$. As $K$ and $K'$ are connected, $K \cup K' \subset U_i$. By part (1), for every $r \in (0, r_i)$ the component $U_r$ containing $K \cup K'$ satisfies $\overline{U_r} \subset X$.
    
Define $C:=\cap_i \overline{U_{r_i}}$. By the same argument as in (1), $C$ is a nonempty, compact, and connected set containing $K \cup K'$. However, since $f$ is continuous, it follows that $f(C) \subset \{y\}$. Therefore, $C$ is contained in a single component of $f^{-1}(y)$, which is only possible if $K = K'$. 
\end{proof}

    We use the following lemma in the next section to prove that our degree function is well-defined on certain subsets of our domain. 
\begin{lemma}
    \label{3.8}
     Let $\varphi:\overline{X} \to \overline{Y}$ be a homeomorphism in $W^{1,n}(X, \mathbb R^n)$. Suppose $f \in W^{1,n}(X, \mathbb R^n)$ is a continuous map on $\overline{X}$ such that $f|_{\partial X} = \varphi|_{\partial X}$. Let $A \subset \mathbb R^n$ such that $\overline{A} \cap \varphi(\partial X) = \emptyset $ and let $C \subset f^{-1}(A)$ be connected. Then $\overline{C} \subset X$. 
\end{lemma}
\begin{proof}
    Suppose for the sake of contradiction that $\overline{C} \cap \partial X \neq \emptyset$. Since $f$ is continuous on $\overline{X}$, we have $\varphi(\overline{C} \cap \partial X) \cap \overline{f(C)} \neq \emptyset$. However, $f(C) \subset A$, so this implies 
    \[\emptyset \neq \varphi(\overline{C} \cap \partial X) \cap \overline{f(C)} \subset \varphi(\partial X) \cap \overline{A} = \emptyset\] which is a contradiction. 
\end{proof}

\section{Continuity up to the Boundary}

We rely on the fact that $X$ is a Lipschitz domain in order to extend $f:X \to \mathbb R^n$ to $\overline{X}$. Then, we study the oscillation of our map to obtain an oscillation estimate. The oscillation of a map $f$ over any subset of the domain is defined:
\begin{definition}
    Let $A \subset \mathbb R^n$ and let $f: A \to \mathbb R^n$. The oscillation of $f$ over a set $A$ is \[\osc_A f := \sup \{|f(x) - f(y)|:x,y \in A\}\]
\end{definition}

\begin{remark}\label{rmk - homeo oscillation}
    For any homeomorphism $\varphi: \overline{X} \to \overline{Y}$, and any $B_r(x) \subset X$ we have $\displaystyle\osc_{B_r(x)} \varphi \le \osc_{\partial B_r(x)} \varphi$. Moreover, this property also holds for any continuous Sobolev map that is weakly monotone~\cite{Hencl}.
\end{remark}
 
 To prove that the extended map $\widetilde{f}: \overline{X} \to \mathbb R^n$ is continuous, we prove an oscillation estimate for the extended map $\widetilde{f}$ and then apply a Sobolev embedding theorem on almost every ball where this estimate holds. 

\begin{prop}\label{7.1}
    Let $\varphi:\overline{X} \to \overline{Y}$ be a given homeomorphism in $W^{1,n}(X, \mathbb R^n)$. Suppose $f \in W^{1,n}(X, \mathbb R^n)$ with $J_f>0$ such that $f - \varphi \in W^{1,n}_0(X, \mathbb R^n)$. Then $f$ extends continuously as a map $f:\overline{X} \to \mathbb R^n$.
\end{prop}
\begin{proof}
    Let $x \in \partial X$ and set $y := \varphi(x)$. Since $Y$ is a Lipschitz domain, there exists some $r_y>0$ and a bi-Lipschitz map $L_y$ that maps $y$ to 0 and $B_{r_y}(y) \cap \partial Y$ onto an $(n-1)$-dimensional hyperplane, $\{(y_1, \dots , y_n) \in B_{r_y}(0) : y_n=0\}$ and $L_y(B_{r_y}(y) \cap Y)= \{y \in B_{r_y}(0) : y_n > 0\}$. Similarly, we can choose some $r_x > 0$ and Lipschitz map $L_x$ so that $L_x(x) = 0$ and $L_x$ maps $B_{r_x}(x) \cap \partial X$ onto an $(n-1)$-dimensional hyperplane, $\{(x_1, \dots , x_n) \in \mathbb R^n : x_n=0\}$. Further, $r_x$ can be chosen so that $\varphi(B_{r_x}(x) \cap \overline{X}) \subset B_{r_y}(y) \cap \overline{Y}$ because of uniform continuity of $\varphi$ on $\overline{X}$.

            We identify $B_{r_x}(x)$ and $B_{r_y}(y)$ with their images under $L_x$ and $L_y$ respectively. Using a standard reflection argument, we can then define 
        \[\widetilde{\varphi} = \begin{cases}
            \varphi & B_{r_x}(0) \cap \{x_n > 0\} \\ (\varphi_1(x_1, \dots , -x_n), \dots, -\varphi_n(x_1, \dots, -x_n)) & B_{r_x}(0) \cap \{x_n \leq 0\}
        \end{cases}\]
        where $\varphi_i$ denotes the coordinate functions of $\varphi$. Using the ACL characterization, $\widetilde{\varphi} \in W^{1,n}(B_{r_x}(0), \mathbb R^n)$. Then define
        \[\widetilde{f} := \begin{cases}
            f &B_{r_x}(0) \cap \{ x_n > 0 \} \\ \widetilde{\varphi} & B_{r_x}(0) \cap \{x_n \leq 0\}
        \end{cases}\]

    Now, $\widetilde{f} - \widetilde{\varphi}=f-\varphi$ on $B_{r_x}(0) \cap \{x_n>0\}$ and 0 otherwise. Since $f-\varphi \in W^{1,n}_0(X, \mathbb R^n)$, we have $\widetilde{f}- \widetilde{\varphi} \in W^{1,n}_0(B_{r_x}(0), \mathbb R^n)$.

    We will show that $\widetilde{f_i}$ for $1 \le i \le n$ is weakly monotone on $B_s(0)$ for a.e. $s \in (0, r_x)$. Let $m<M$ be real numbers and assume that $v:=(\widetilde{f_i}-M)^+ - (m-\widetilde{f_i})^+ \in W^{1,n}_0(B_{s}(x), \mathbb R^n)$. Then 
    \[\triangledown v = \begin{cases}
        0 & m \le \widetilde{f_i} \le M \\ \triangledown \widetilde{f_i} & \text{otherwise}
    \end{cases}\]
Let $E$ denote the set where $\triangledown v = \triangledown \widetilde{f_i}$. Since $\widetilde{\varphi}$ is a homeomorphism, the maximum/minimum of $\widetilde{\varphi_i}$ is attained on $\partial B_s(0)$, so $\triangledown v =0$ when $\widetilde{f_i} = \widetilde{\varphi_i}$. Consequently, $E \subset B_{s}(0) \cap \{x_n>0\} $, so $J_{\widetilde{f}} = J_f >0$ in $E$. Then
 \begin{align*}
    \int_{B_s(0)} \frac{|\triangledown v|^n}{|D\widetilde{f}(x)|^n }J_{\widetilde{f}}(x) dx &= \int_E  \frac{|\triangledown \widetilde{f_i}|^n}{|D\widetilde{f}(x)|^n }J_f(x) dx \\ &\le \int_E J_f(x) dx 
\end{align*}

Now, \begin{align*}
    \int_E J(f_1, \dots, f_n)(x) dx &= \int_E J(f_1, f_2, \dots , f_{i-1}, v, f_{i+1}, \dots, f_n) dx \\ &= \int_{B_s(0)} J(f_1, f_2, \dots , f_{i-1}, v, f_{i+1}, \dots, f_n) dx
\end{align*} because of how $v$ is defined. Further, $v \in W^{1,n}_0(B_s(0))$, so by Stokes' Theorem, $\displaystyle\int_{B_s(0)} J(f_1, f_2, \dots , f_{i-1}, v, f_{i+1}, \dots, f_n) dx =0$.

Thus, we've shown $\triangledown v=0$ almost everywhere in $B_s(0)$, which implies $v$ is constant, but since $v \in W^{1,n}_0(B_s(0))$, $v \equiv 0$ in $B_s(0)$. Hence $m \le \widetilde{f_i} \le M$ almost everywhere in $B_s(0)$.

For each $1 \le i \le n$, $\widetilde{f_i}$ is a real-valued weakly monotone Sobolev mapping, all of which satisfy the following oscillation estimate $\displaystyle \osc_{B_{r_x/2}(0)}\widetilde{f_i} \le \osc_{\partial B_{r_x/2}(0)} \widetilde{f_i}$.
Then, since norms in a finite-dimensional vector space are equivalent,
        \begin{align*}
            \osc_{B_{r_x/2}(0)}\widetilde{f} &\le \sum_{i=1}^n \osc_{B_{r_x/2}(0)} \widetilde{f}_i \\&\le \sum_{i=1}^n \osc_{\partial B_{r_x/2}(0)} \widetilde{f}_i \\&\le \sqrt{n} \osc_{\partial B_{r_x/2}(0)}\widetilde{f}.
        \end{align*}
Lastly, we use the Sobolev embedding theorem on spheres~\cite[Theorem A.18]{Hencl} to demonstrate continuity. Indeed, let $r \in (0, r_x/4)$. For a.e. $r < s < r_x/2$ we have 
\begin{align*}
        \osc_{B_r(x)} \widetilde{f} &\le \osc_{B_s(x)} \widetilde{f} \\ &\le \sqrt{n} \osc_{\partial B_s(x)} \widetilde{f} \\ &\le C(n) s^{1/n} \left( \int_{\partial B_s(x)} |D\widetilde{f}|^n \right)^{1/n}
    \end{align*}
    which implies that 
        \[\frac{[\osc_{B_r(x)}\widetilde{f}]^n}{s} \le C(n)\int_{\partial B_s(x)} |D\widetilde{f}|^n\]

        Since the above holds for a.e. $s \in (r, r_x/2)$ we can integrate as follows
        \[\int_r^{r_x/2} \frac{[\osc_{B_r(x)}\widetilde{f}]^n}{s} ds  \le C(n) \int_r^{r_x/2} \int_{\partial B_s(x)} |D\widetilde{f}|^n ds\]

        Observe that the right hand side is equivalent to integrating $|D\widetilde{f}|^n$ over the annulus $B_{r_x/2} \setminus B_r$, so we have that 
        \[\int_r^{r_x/2} \frac{[\osc_{B_r(x)}\widetilde{f}]^n}{s} ds  \le C(n) \int_{ B_{r_x/2}(x)} |D\widetilde{f}|^n ds\]
Thus for $r \in (0, r_x/2)$ we have \begin{equation}\label{eqn-Kosc}[\osc_{B_r(x)} \widetilde{f}]^n \le \frac{C}{\log \left( \frac{r_x}{2r}\right)}\int_{B_{r_x/2}(x)} |D\widetilde{f}|^n\end{equation} and, in particular, $\widetilde{f}$ is continuous on $B_r(x)$ for a.e. $r \in (0, r_x/2)$.
\end{proof}

\section{Properties of the extended map}\label{sec-properties}
In this section, we exclusively work with the extended map in Proposition~\ref{7.1}, which we still refer to as $f$. We prove that our extended map $f: \overline{X} \to \mathbb R^n$ inherits some properties of $\varphi$, such as surjectivity. However, as we will see in Example~\ref{mtexample}, maps defined in Theorem~\ref{thm-main} may not have injectivity in the classical sense. Instead, such maps will have global invertibility in the sense of Ball~\cite{ball} and monotonicity in the sense of Morrey~\cite{morrey}.

A primary tool in this section is degree theory. Since we've shown that $f$ extends continuously to the boundary, we can use a well-known property of the degree function that if $f=\varphi$ on $\partial X$ then 
\begin{equation}\label{eqn-deg boundary}
    \deg(f, y, X) = \deg(\varphi, y, X) 
\end{equation}
for all $y \in \mathbb R^n \setminus \partial Y$

Another significant property of the degree function is the ability to compute the degree over a disconnected set by computing the degrees over each of its components: 

\begin{definition}\label{def-degreesum} \cite[Theorem 2.7]{Fonseca}
    Let $U \subset \mathbb R^n$ be an open subset and let $f: \overline{U} \to \mathbb R^n$ be continuous. Let $U = \cup_{i \in \mathbb N}U_i$ be a union of open and mutually disjoint $U_i$. Then 
    \[\deg(f, y, U) := \sum_i \deg(f, y, U_i)\] for any $y \notin f(\partial U)$
\end{definition}

First, we prove a useful application of Theorem~\ref{Fonseca}. This proof does not require $f$ to be continuous on $\overline{X}$, so we've stated it in full generality.
\begin{lemma}\label{5.1} 
Let $\varphi: \overline{X} \to \overline{Y}$ be a homeomorphism in $W^{1,n}(X, Y)$. Suppose $f \in W^{1,n}(X, \mathbb R^n)$ with $J_f > 0 $ a.e. and $f - \varphi \in W^{1,n}_0(X, \mathbb R^n)$. Let $V \subset \mathbb R^n$ be an open, bounded, connected set such that $\dist(\partial Y,\overline{V}) > 0$ and $|\partial V| =0$. Let $U $ be a connected component of $f^{-1}(V)$. Then $\overline{U} \subset X$ and for every $y \in V$, 
\begin{equation}\label{eq-degree}
    \deg (f, y, U) = \frac{1}{|V|} \int_U J_f(x) dx 
    \end{equation}
\end{lemma}
\begin{proof}
 Since $\dist(\partial Y,\overline{V}) > 0$, we have that $\partial Y \cap \overline{V}  = \emptyset$, so by Lemma~\ref{3.8}, $\overline{U} \subset X$.

First, we show $\deg(f, y, U)$ is well-defined for all $y \in V$:
As a connected component of $f^{-1}(V)$, U is both closed and open in $f^{-1}(V)$, so $\partial U \cap f^{-1}(V)=\emptyset$. 

This implies $f(\partial U) \subset \mathbb R^n \setminus V$. However, since $f$ is continuous on $X$ and $\overline{U} \subset X$, we have that $f(\partial U) \subset \overline{f(U)} \subset \overline{V}$. Thus, $f(\partial U) \subset \partial V$ so $\deg (f, y, U)$ is well-defined for all $y\in V$.

Now, $|\partial U|=0$ because $|\partial V|=0$ and $f$ has the Lusin $N^{-1}$ condition by Remark~\ref{conditionN}. By Theorem~\ref{Fonseca}
\[\int_U J_f(x) dx = \int_{\mathbb R^n} \deg(f, y', U) dy' \]
 which can be simplified to

\[ \int_U J_f(x) dx  = \int_V \deg(f, y', U) dy'   \] because $f(U) \subset V$, so $\deg(f, y, U) =0$ for any $y \in \mathbb R^n \setminus \overline{V}$. The degree map is constant on connected sets~\cite[Theorem 2.3(3)]{Fonseca}, so for all $y \in V$, 
\[\int_U J_f(x) dx = |V|\deg(f, y, U)  \]
\end{proof}
\begin{remark}\label{degreermk}
    The conditions $\varphi -f \in W^{1,n}_0(X, \mathbb R^n)$ and $\dist(\partial Y,\overline{V}) > 0$ are only needed so that $\overline{U} \subset X$. If apriori we have $\overline{U} \subset X$ then the same proof applies if $f$ is a continuous map on $X$ with Lusin $N^{-1}$ condition and $V$ is an open, bounded, connected set with $|\partial V|=0$. 
\end{remark}
Using the multiplicity function $N(y, f, X) := \#\{x \in f^{-1}(y) \cap X\}$ and its relationship with the degree function, we prove surjectivity and global invertibility of $f$.
\begin{lemma}\label{onto}
        Let $\varphi: \overline{X} \to \overline{Y}$ be a homeomorphism in $W^{1,n}(X,Y)$. Let $f \in W^{1,n}(X, \mathbb R^n)$ such that $f $ is continuous on $\overline{X}$, $f=\varphi$ on $\partial X$ and $J_f > 0$ a.e. Then $f$ maps $\overline{X}$ onto $\overline{Y}$.
\end{lemma}
\begin{proof}
    By (\ref{eqn-deg boundary}), $\deg(f, y, X) =1$ for all $y \in Y$, so by~\cite[Theorem 2.4]{Fonseca}, $Y \subset f(X)$. 

     To show the opposite containment, let $y \in \mathbb R^n \setminus \overline{Y}$, and assume, for the sake of contradiction, that there is some $x \in X$ such that $f(x) = y$.

 By Lemma~\ref{lem-component}, there exists some $r> 0$ such that $f^{-1}(B_r(y)) \subset X$, and, shrinking if necessary, we can choose $r$ such that $B_r(y) \cap \overline{Y} = \emptyset$.
 
 Using local invertibility of $W^{1,n}$ mappings as in~\cite[Theorem 6.1]{Fonseca} there exists an open set $U$ around $y$ such that $N(y', f, X) \ge 1$ for a.e. $y' \in U$. Consider the open set $B_r(y) \cap U$. Combining~\cite[Theorem 5.27]{Fonseca} and~\cite[Theorem 5.34]{Fonseca} we have 
\[\int_{B_r(y) \cap U} N(y', f, X) dy'= \int_{B_r(y) \cap U} \deg(f, y', X) dy' =0\] where the last equality follows from (~\ref{eqn-deg boundary}) because $y' \in f(X) \setminus Y$. However, since $N(y', f, X) \ge 1$ for almost every $y'\in B_r(y) \cap U$, this implies $|B_r(y) \cap U| = 0$, which is a contradiction. Hence, $f$ maps $\overline{X}$ onto $\overline{Y}$.\\ 
\end{proof}

\begin{lemma}\label{lem-globallyinvertible}
    Let $\varphi: \overline{X} \to \overline{Y}$ be a homeomorphism in $W^{1,n}(X,Y)$. Let $f \in W^{1,n}(X, \mathbb R^n)$ such that $f $ is continuous on $\overline{X}$, $f=\varphi$ on $\partial X$ and $J_f > 0$ a.e. Then $f$ is globally invertible.
\end{lemma}

\begin{proof}
We prove global invertibility by showing (\ref{globalinvertibility}). Previously, we saw that (\ref{eqn-deg boundary}) implies $Y \subset f(X)$. Thus, for any $y \in Y$, we have that $N(y, f, X) \ge 1$. For the sake of contradiction, assume that the set $Y_0 := \{y \in Y : N(y, f, X) > 1\}$ has positive measure. Then by~\cite[Theorem A.35]{Hencl}
    \[\int_X J_f(x) dx = \int_{Y\setminus Y_0} N(y, f, X) dy + \int_{Y_0} N(y, f, X) dy > |Y|\] 

    However, because $f$ and $\varphi$ are $W^{1,n}$ mappings that agree  on $\partial X$, it follows from Theorem~\ref{Fonseca} that
    \[\int_X J_f(x) dx = \int_X J_{\varphi}(x) dx = |Y|\]
    which is a contradiction. Thus $|Y_0|=0$, so (\ref{globalinvertibility}) holds.
\end{proof}

\begin{proof}[Proof of Corollary~\ref{cor-injective}]
   It follows from Lemma~\ref{lem-globallyinvertible} and $f$ having the Lusin $N^{-1}$ condition that $X_f := X \setminus f^{-1}(Y_0)$ has full measure in $X$. 
\end{proof}
To prove monotonicity, we will consider interior points and boundary points separately. 
\begin{lemma}
        Let $\varphi: \overline{X} \to \overline{Y}$ be a homeomorphism in $W^{1,n}(X,Y)$. Let $f \in W^{1,n}(X, \mathbb R^n)$ such that $f $ is continuous on $\overline{X}$, $f=\varphi$ on $\partial X$ and $J_f > 0$ a.e. Then $f^{-1}(y)$ is connected for all $y \in Y$.
\end{lemma}

\begin{proof}
    Let $y_0 \in Y$ and suppose that $f^{-1}(y_0)$ is disconnected with components $K_1, K_2,\dots$. Since the $K_i$ are open and disjoint, there are at most countably many of them. Since $y_0$ is an interior point, $\dist(y_0, f(\partial X)) > 0$ and for sufficiently small $r$, $\overline{f^{-1}(B_r(y_0))} \subset X$. Observe that $f(\partial f^{-1}(B_r(y_0))) \subset \partial f(f^{-1}(B_r(y_0))) = \partial B_r(y_0)$ by Lemma~\ref{onto}. Thus, $y_0 \notin f(\partial(f^{-1}(B_r(y_0)))$ and the degree is well defined on $f^{-1}(B_r(y_0))$ and its components of which there are at most countably many. Let $U_i$ denote the components of $f^{-1}(B_r(y_0))$ that contain $K_i$.
    We can now use Definition~\ref{def-degreesum} to write
   \[\deg(f, y_0, f^{-1}(B_r(y_0))) = \deg (f, y_0, U_r)  + \deg (f, y_0, V_r) \] 

    By (\ref{eqn-deg boundary}) $\deg(f, y_0, X) = 1$ because $y_0 \in Y$. Observe that $X \setminus  f^{-1}(B_r(y_0))$ contains no preimages of $y_0$ under $f$. By~\cite[Theorem 2.1]{Fonseca} and Definition~\ref{def-degreesum}, $\deg (f, y_0, f^{-1}(B_r(y_0))) = \deg (f, y_0, X)$. Thus, $\deg (f, y_0, f^{-1}(B_r(y_0))) = 1 $.  
    
By Lemma~\ref{5.1} and (\ref{eq-degree})
\[\int_{U_r} J_f(x) dx = |B_r(y_0)| \deg(f, y_0, U_r)\] Thus, $\deg(f, y_0, U_r) > 0$ and similarly, $\deg(f, y_0, V_r) >0$. However, the degree function is integer-valued, so $\deg (f, y_0, f^{-1}(B_r(y_0))) \ge 2 $ which is a contradiction. 
\end{proof}

To show the boundary case, we require an extra lemma. 
   \begin{lemma}\label{6.4} Let $\varphi:\overline{X} \to \overline{Y}$ be a homeomorphism in $W^{1,n}(X, \mathbb R^n)$. Let $f \in W^{1,n}(X, \mathbb R^n)$ such that $f$ is continuous on $\overline{X}$, $f_{\partial X}=\varphi|_{\partial X}$, and $J_f > 0$ a.e. If $y\in \partial Y$, $x_0=\varphi^{-1}(y)$, and $K$ is a connected component of $f^{-1}(y)$ then $x_0  \in \overline{K}$. 
\end{lemma}
\begin{proof} 
We proceed by breaking this proof into cases. \\
    Case 1: Assume that $\overline{K} \cap \partial X \neq \emptyset$. Then since $f|_{\partial X}=\varphi|_{\partial X}$, we have $\varphi(\overline{K}) \cap f(\overline{K}) \neq \emptyset$. As $K$ is a component of $f^{-1}(y)$ then $f(\overline{K})= \{y\}$. Since $f(\overline{K})$ contains one element and intersects the set  $\varphi(\overline{K})$, it must be the case that $y \in \varphi(\overline{K})$, so $x_0 \in \overline{K}$.

    Case 2: Assume $\overline{K} \cap \partial X = \emptyset$ and, by way of contradiction, that $x_0 \notin \overline{K}$. By Lemma~\ref{lem-component}, there exists some $r> 0$ and component $U$ of $ f^{-1}( B_r(y))$ such that $K \subset U$ and $\overline{U} \subset X$. Using Remark~\ref{degreermk},    
    \[\frac{1}{|B_r(y)|}\int_U J_f(x) dx = \deg(f, y, U)\]
    which implies $\deg(f, y, U) > 0$, so by~\cite[Theorem 2.1]{Fonseca} 
    $B_r(y) \subset f(U)$. Now, since $y$ is on the boundary of a Lipschitz domain, $B_r(y) \cap Y \neq \emptyset$ and $B_r(y) \cap \mathbb R^n \setminus Y \neq \emptyset$. However, $f(X) \subset \overline{Y}$ by Lemma~\ref{onto}, so it must be the case that $B_r(y) \cap \mathbb R^n \setminus Y = \emptyset$, which is a contradiction.
\end{proof}

\begin{cor}\label{6.5}
Let $\varphi:\overline{X} \to \overline{Y}$ be a homeomorphism in $W^{1,n}(X, \mathbb R^n)$. Suppose $f \in W^{1,n}(X, \mathbb R^n)$ is a continuous map on $\overline{X}$ such that $J_f > 0$ a.e. and $f|_{\partial X} = \varphi|_{\partial X}$. Then for every boundary point $y \in \partial Y$, the set $f^{-1}(y)$ is connected.
\end{cor}
\begin{proof}
    
Suppose that the open sets $(U,V)$ form a separation of $f^{-1}(y)$. Without loss of generality, say $x_0 \in U$. 

   Let $K$ be a connected component of $f^{-1}(y)$. By Lemma~\ref{6.4}, $x_0 \in \overline{K}$, and since $U$ is an open set containing $x_0$, $U \cap K \neq \emptyset$. As $K$ is connected, $V \cap K = \emptyset$. Since $K$ was arbitrary, $V = \emptyset$. Therefore, $f^{-1}(y) \cup (f|_{\partial X})^{-1}\{y\}$ is connected. 
\end{proof}

\section{Minimizing a Neohookean Energy Functional}\label{app}
As an application of our results, we use the Direct Method of Calculus of Variations to solve a minimization problem for a Neohookean energy functional: 
\[E^p_{\alpha}[f] := \int_B |Df|^p + \frac{1}{J_f^{\alpha}}\qquad \alpha >0\]

The term Neohookean refers to the fact that $E^p_{\alpha}[f]$ approaches $\infty$ as $J_f \to 0^+$. In physical applications, this condition represents how an infinite amount of work is required to compress a finite volume of material into zero volume.

Given a homeomorphism $\varphi: \overline{X} \to \overline{Y}$, we can define an admissible class of deformations \[\mathcal{A}^p_{\varphi} = \{f \in W^{1,p}(X, \mathbb R^n) : J_f > 0  \text{ a.e.}, f-\varphi \in W^{1,p}_0(X, \mathbb R^n)\}\]
 When $p=n$, Theorem~\ref{thm-main} shows that such mappings are continuous on $\overline{X}$, surjective, and monotone. Global invertibility is also shown but will not be featured in this section.

The main result of this section is that there exists some mapping in $\mathcal{A}^p_{\varphi}$ which minimizes $E^p_{\alpha}[f]$ subject to $\mathcal{A}^p_{\varphi}$:

\begin{theorem}\label{energyminthm}
    Let $p \ge n$. There exists an $f_0 \in \mathcal{A}^p_{\varphi}$ such that \[\inf_{f \in \mathcal{A}^p_{\varphi}} E[f] = E[f_0] \]
\end{theorem}

\begin{remark}
    It follows from \cite{ball} and Theorem~\ref{thm-main} that when $p \ge n$, minimizing $E^p_{\alpha}[f]$ subject to $\mathcal{A}^p_{\varphi}$ is the same as minimizing subject to $\{f: \overline{X} \to \overline{Y} : f \in W^{1,p}(X, \mathbb R^n), J_f >0 $ a.e., monotone, and $f=\varphi$ on $\partial X\}$.
\end{remark}
\begin{lemma}\label{wlsc} 
    For sequences in $\mathcal{A}^n_{\varphi}$ the energy functional is weakly lower semi-continuous, i.e. for $f_k \rightharpoonup f$
    \[\liminf_{k \to \infty} E[f_k] \ge E[f] \]
\end{lemma}

\begin{proof}
    Assume $\{f_k\} \subset \mathcal{A}^n_{\varphi}$  converges weakly in $W^{1,n}(X, \mathbb R^n)$ to some $f$. 

    As the $L^n$ norm is lower semi-continuous with respect to the weak topology 
    \[\|Df\|_n \le \liminf_{k \to \infty} \|Df_k\|_n\] so to show the claim, it suffices to show 
    \[\liminf_{k \to \infty} \| (Jf_k)^{-\alpha}\| \ge \| (Jf)^{-\alpha}\|\]

    Now, since $t^{-\alpha}$ is convex for $\alpha > 0$, for all $\varepsilon >0$ we have that
    \begin{equation}\label{eq:1} \frac{1}{(Jf_k + \varepsilon)^{\alpha}} -\frac{1}{(Jf + \varepsilon)^{\alpha}} \ge \frac{-\alpha}{(J_f + \varepsilon)^{\alpha + 1} } (Jf_k - Jf)\end{equation} Additionally, $J_f >0$ a.e. and $\varepsilon >0$ implies that $\displaystyle\frac{-\alpha}{(J_f + \varepsilon)^{\alpha + 1} } \in L^{\infty}(X)$.
    
    By~\cite[Theorem 8.4.2]{Tadeuszbook}, since $J_{f_k} >0$ a.e. the sequence of Jacobians $\{J_{f_k}\}$ are weakly convergent in $L^1_{\text{loc}}(X, \mathbb R^n)$, so the following limit holds
    \[\lim_{k \to\infty} \int_X\frac{-\alpha}{(J_f + \varepsilon)^{\alpha + 1} }(Jf_k - Jf) = 0\] 

    Combining the above with Equation~\ref{eq:1} we see
    \[\liminf_{k \to \infty} \int_X\frac{1}{(Jf_k)^{\alpha}} \ge \int_X\frac{1}{(Jf + \varepsilon)^{\alpha}}\]

    Since this is true for all $\varepsilon > 0$, select a decreasing sequence of $\varepsilon_n \searrow 0$ so that by Monotone Convergence Theorem 
    \[\liminf_{k \to \infty} \int_X\frac{1}{(Jf_k)^{\alpha}} \ge \int_X\frac{1}{(Jf )^{\alpha}}\]

    Hence 
    \begin{align*}
        \liminf_{k \to \infty} E[f_k] &= \liminf_{k \to \infty} \int_X \left( |Df_k|^n + \frac{1}{(Jf_k)^{\alpha}} \right) \\ &\ge \int_X \left(|Df|^n + \frac{1}{J_f^{\alpha} } \right) \\ &=E[f]
    \end{align*}
\end{proof}

\begin{proof}[Proof of Theorem~\ref{energyminthm}]
The $p > n$ case is proven in~\cite{ball}. Let $\{f_k\} \subset \mathcal{A}^n_{\varphi}$ be a minimizing sequence. By lemma~\ref{wlsc}, $E[f] < \infty$ so $J_f > 0$ a.e. 
The remaining properties follow from Lemma~\ref{wlsc}, and the following chain of inequalities completes the proof.
\[\inf_{f \in \mathcal{A}^p_{\varphi}} E[f] \le E[f_0] \le \liminf_{k_j \to \infty} E[f_{k_j}] =\inf_{f \in \mathcal{A}^p_{\varphi}} E[f] \]
\end{proof}
\section{Examples}\label{examples}
 \begin{example}\label{ex-ctyonto}
    Let $\mathbb D$ be the unit disk in $\mathbb C$ and define a sequence of disjoint balls $B_k$ ($k=1, 2, \dots) $ with centers $\omega_k : = 1 -\frac{1}{k}$ and radii $2s_k$ where $s_k := 10^{-k}$. On each $B_k$, define 
    \[G_k(\omega_k + re^{i\theta}) = \begin{cases} \omega_k + 2(r-s_k)e^{i\theta} & s_k \le r \le 2s_k \\ \omega_k + \frac{\log(s_k/r)}{k^4} & s_ke^{-k^{4}} < r < s_k \\ \omega_k + 2 & r \le s_ke^{-k^{4}}
    \end{cases}\]

    We claim the map $G: \overline{\mathbb D} \to \overline{\mathbb C}$ 

    \[G(z) = \begin{cases} G_k(z) & z \in B_k, k=1, 2, \dots \\ z & \text{otherwise}
    \end{cases}\]
is not surjective nor can it be extended continuously to $S^1(0)$ when we allow the Jacobian to be 0 on a set of positive measure. 
First, observe that 
\begin{align*}
    \|G\|_{W^{1,2}(\mathbb D, \mathbb C)} &\le C + \sum_{k=1}^\infty \|DG_k(z)\|_{L^2(B_k)} \\ &\le C + \sum_{k=1}^{\infty } \sqrt{C_1 10^{-2k} + \frac{8\pi}{k^4}} \\ &\le C+C_2 \sum_{k=1}^\infty \left(10^{-k} + \frac{1}{k^2}\right) < \infty
\end{align*}
so $G \in W^{1,2}(\mathbb D, \mathbb C)$ and moreover, since $\displaystyle \lim_{x \to \omega} G(z) = \omega $ for all $\omega \in \partial \mathbb D \setminus \{1\}$, the trace of $G$ on the boundary must also equal the identity mapping. Thus, $G - \text{id} \in W^{1,2}_0 (\mathbb D, \mathbb C)$.

    The mapping is clearly continuous on $\mathbb D$ but since $G$ maps the centers $w_k$ of each $B_k$ to $w_k + 2$, which has real part larger than 2, it fails to be continuous up to the boundary. Additionally, since $G(\mathbb D) \setminus \overline{ \mathbb D} \neq \emptyset$, $G$ is not a surjective map. 
\end{example}
\begin{example}\label{ex-ctybdry}
    Let $\mathbb D$ again be the unit disk in $\mathbb C$ and define $G:\overline{\mathbb D} \to \overline{\mathbb D} $ as 
\[G(re^{i\theta}) = \begin{cases} 2(r-1/2)e^{i\theta} & 1/2 \le r \le 1 \\ 1-2r & r < 1/2
\end{cases}\]
Observe that $G$ is the identity on the boundary and $J_G \ge 0$ a.e. However, for any point $z \in (0,1) \times \{0\}$, $G^{-1}(z)$ is the union of a single point in $ \mathbb D \setminus \overline{B(0,1/2)} $ and a circle of radius $r_0 \in (0, 1/2)$, causing monotonicity to fail.
\end{example} 
\begin{example}\label{mtexample}
    Consider the rectangle $R = \{(x,y) \in \mathbb R^2 : |x| \le 1, |y| \le 2\}$. Define $H: R \to \mathbb R^2$ 
    \[ H(x_1, x_2) = \begin{cases}
        (x_1, |x_1|x_2) & |x_2|\le 1 \\ (x_1, [2(|x_2|-1)+(2-|x_2|)|x_1|]\sgn(x_2) & 1 < |x_2| \le 2
    \end{cases}\]
    Despite having an almost everywhere positive Jacobian and being the identity on the boundary $\partial R$, $H$ is not injective.
\end{example}

\end{document}